\documentclass[submission]{FPSAC2018}

\usepackage{pstricks,pst-node,graphicx}

\newtheorem{theorem}{Theorem}

\newtheorem{proposition}[theorem]{Proposition}
\newtheorem{corollary}[theorem]{Corollary}

\theoremstyle{definition}
\newtheorem{definition}[theorem]{Definition}

\theoremstyle{remark}

\numberwithin{equation}{section}

\newcommand{\st}[1]{\ensuremath{#1}^{\prime}}

\newcommand{\wt}{\ensuremath\mathrm{wt}}
\newcommand{\SSYT}{\ensuremath\mathrm{SSYT}}
\newcommand{\SSHT}{\ensuremath\mathrm{SSHT}}
\newcommand{\SSST}{\ensuremath\mathrm{SSST}}

\newcommand{\Yam}{\ensuremath\mathrm{Yam}}
\newcommand{\sym}{\ensuremath\mathrm{sym}}
\newcommand{\rect}{\ensuremath\mathrm{rect}}

\newlength\cellsize 

\savebox2{%
\begin{picture}(12,12)
\put(0,0){\line(1,0){12}}
\put(0,0){\line(0,1){12}}
\put(12,0){\line(0,1){12}}
\put(0,12){\line(1,0){12}}
\end{picture}}

\newcommand\cellify[1]{\def\thearg{#1}\def\nothing{}%
\ifx\thearg\nothing\vrule width0pt height\cellsize depth0pt%
  \else\hbox to 0pt{\usebox2\hss}\fi%
  \vbox to 12\unitlength{\vss\hbox to 12\unitlength{\hss$#1$\hss}\vss}}

\newcommand\tableau[1]{\vtop{\let\\=\cr
\setlength\baselineskip{-12000pt}
\setlength\lineskiplimit{12000pt}
\setlength\lineskip{0pt}
\halign{&\cellify{##}\cr#1\crcr}}}

\savebox4{
\begin{picture}(12,12)
\put(6,6){\circle{12}}
\end{picture}}

\newcommand{\cir}[1]{\def\thearg{#1}\def\nothing{}%
\ifx\thearg\nothing\vrule width0pt height\cellsize depth0pt%
  \else\hbox to 0pt{\usebox4\hss}\fi%
  \vbox to 12\unitlength{\vss\hbox to 12\unitlength{\hss$#1$\hss}\vss}}

\definecolor{boxgray}{gray}{.8}
\newcommand{\cb}{\color{boxgray}\rule{1\cellsize}{1\cellsize}\hspace{-\cellsize}\usebox2}

\newcommand{\tb}{%
  \psset{unit=\cellsize}
  \begin{pspicture}(1,1)
    \psset{linewidth=0.1ex}
    \pscustom[fillstyle=solid,fillcolor=boxgray]{
      \pspolygon(0,0)(1,1)(0,1)}
  \end{pspicture}}

\usepackage{lipsum}

\title{Crystal graphs for shifted tableaux}

\author{Sami Assaf\thanks{\href{mailto:shassaf@usc.edu}{shassaf@usc.edu} Work supported by a grant from the Simons Foundation (Award 524477, S.A.).}\addressmark{1} \and Ezgi Kantarc\i~O\u{g}uz\thanks{\href{mailto:kantarci@usc.edu}{kantarci@usc.edu}}\addressmark{1}}

\address{\addressmark{1}Department of Mathematics, University of Southern California, 3620 South Vermont Avenue, Los Angeles, CA 90089-2532, U.S.A.}

\received{\today}

\abstract{We define crystal operators on semistandard shifted tableaux, giving a new proof that Schur $P$-functions are Schur positive. We define a queer crystal operator to construct a connected queer crystal on semistandard shifted tableaux of a given shape, providing a new proof that products of Schur $P$-functions are Schur $P$-positive. We also give a rectification map from shifted tableaux to Young tableaux that commutes with the crystal operators and provides a dual algorithm to shifted insertion.}


\usepackage[backend=bibtex]{biblatex}
\begin{document}

\maketitle

%
\section{Introduction}
%
\label{sec:introduction}

Schur polynomials are ubiquitous throughout mathematics, arising as irreducible characters for polynomial representations of the general linear group, Frobenius characters for irreducible representations of the symmetric group, and polynomial representatives for the cohomology classes of Schubert cycles in Grassmannians. Combinatorially, Schur polynomials are the generating polynomials for semistandard Young tableaux.

Schur $P$-polynomials arise as characters of tensor representations of the queer Lie superalgebra \cite{Ser84}, characters of projective representations of the symmetric group \cite{Ste89}, and representatives for cohomology classes dual to Schubert cycles in isotropic Grassmannians \cite{Pra91}. They enjoy many properties parallel to Schur polynomials. Combinatorially, Schur $P$-polynomials are the generating polynomials for semistandard shifted tableaux.

Stanley conjectured that Schur $P$-polynomials are Schur positive, and this follows from Sagan's shifted insertion \cite{Sag87} independently developed by Worley \cite{Wor84}. Assaf \cite{Ass18} gave another proof using the machinery of dual equivalence graphs \cite{Ass15}. In this extended abstract, we present a new proof using the machinery of crystal graphs \cite{Kas91} and Stembridge's local characterizations thereof \cite{Ste03}, and we provide a dual construction to shifted insertion that \emph{rectifies} a shifted tableau to a Young tableau.

Kashiwara \cite{Kas91} introduced crystal bases to study representations of quantized universal enveloping algebras. A \emph{crystal graph} is a directed, colored graph with vertex set given by the crystal basis and directed edges given by deformations of the Chevalley generators. Crystal graphs encode important information, for example tensor decompositions, for the corresponding representations. The crystal basis for the general linear group can be indexed by semistandard Young tableaux, and there is an explicit combinatorial construction of the crystal graph on tableaux \cite{KN94,Lit95}. Recently, Grantcharov, Jung, Kang, Kashiwara, and Kim \cite{GJKKK14} developed crystal bases for the quantum queer superalgebra and gave an explicit construction of the queer crystal on semistandard decomposition tableaux \cite{Ser10}, another combinatorial model for Schur $P$-polynomials. In this extended abstract, we give an explicit construction of the queer crystal on semistandard shifted tableaux and use this to give a local characterization of queer crystals. 

In Section~\ref{sec:young}, we review Young tableaux and crystal graphs. In Section~\ref{sec:shifted}, we review shifted tableaux and define our crystal operators on semistandard shifted tableaux. In Section~\ref{sec:queer}, we extend our crystal on semistandard shifted tableaux to a queer crystal graph. In Section~\ref{sec:rectify}, we define a \emph{shifted rectification} map from shifted tableaux to Young tableaux that commutes with the crystal operators. For further details, see \cite{AO}.

%
\section{Crystal graph on semistandard Young tableaux}
%
\label{sec:young}

A \emph{partition} $\lambda$ is a weakly decreasing sequence of positive integers, $\lambda = (\lambda_1,\lambda_2, \ldots, \lambda_{\ell})$, where $\lambda_1 \geq \lambda_2 \geq \cdots \geq \lambda_{\ell} > 0$. We identify a partition $\lambda$ with its \emph{Young diagram}, the collection of left-justified cells with $\lambda_i$ cells in row $i$ indexed from the bottom. 

\begin{figure}[ht]
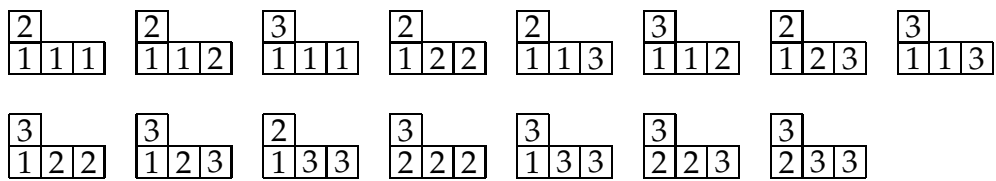

  \begin{center}
    \begin{displaymath}
      \begin{array}{c@{\hskip\cellsize}c@{\hskip\cellsize}c@{\hskip\cellsize}c@{\hskip\cellsize}c@{\hskip\cellsize}c@{\hskip\cellsize}c@{\hskip\cellsize}c}
        \tableau{2 \\ 1 & 1 & 1} & \tableau{2 \\ 1 & 1 & 2} & \tableau{3 \\ 1 & 1 & 1} & \tableau{2 \\ 1 & 2 & 2} & 
        \tableau{2 \\ 1 & 1 & 3} & \tableau{3 \\ 1 & 1 & 2} & \tableau{2 \\ 1 & 2 & 3} & \tableau{3 \\ 1 & 1 & 3} \\ \\
        \tableau{3 \\ 1 & 2 & 2} & \tableau{3 \\ 1 & 2 & 3} & \tableau{2 \\ 1 & 3 & 3} & \tableau{3 \\ 2 & 2 & 2} & 
        \tableau{3 \\ 1 & 3 & 3} & \tableau{3 \\ 2 & 2 & 3} & \tableau{3 \\ 2 & 3 & 3} & 
      \end{array}
    \end{displaymath}
    \caption{\label{fig:Y31}The semistandard Young tableaux of shape $(3,1)$ with entries $\{1,2,3\}$.}
  \end{center}
\end{figure}

The \emph{semistandard Young tableaux} of shape $\lambda$, denoted by $\SSYT(\lambda)$, are fillings of the Young diagram for $\lambda$ with positive integers such that entries weakly increase along rows and columns and each column has at most one entry $i$. For example, see Figure~\ref{fig:Y31}. 

\begin{definition}
  The \emph{Schur polynomial} indexed by the partition $\lambda$ is given by
  \begin{equation}
    s_{\lambda} (x_1,\ldots,x_n) = \sum_{T \in \SSYT_n(\lambda)} x_1^{\wt(T)_1} \cdots x_n^{\wt(T)_n}.
    \label{e:schur_P}
  \end{equation}
  where $\wt(T)$ is the weak composition whose $i$th part is the number of entries $i$ in $T$.
  \label{def:schurP}
\end{definition}
  
For example, from Figure~\ref{fig:S31}, we compute
\begin{eqnarray*}
  s_{(3,1)}(x_1,x_2,x_3) & = & x_1^3 x_2 + x_1^3 x_3 + x_1^2 x_2^2 + 2 x_1^2 x_2 x_3 + x_1^2 x_3^2 + x_1 x_2^3 + 2 x_1 x_2^2 x_3 \\
  & & + 2 x_1 x_2 x_3^2 + x_1 x_3^3 + x_2^3 x_3 + x_2^2 x_3^2 + x_2 x_3^3.
\end{eqnarray*}

For a word $w$ of length $k$, a positive integer $r \leqslant k$, and a positive integer $i$, define
\begin{equation}
  m_i(w,r) = \wt(w_{1} w_{2} \cdots w_{r})_{i} - \wt(w_{1} w_{2} \cdots w_{r})_{i+1}, \hspace{2em}
  m_i(w) = \max_r(m_i(w,r)).
  \label{e:r-index}
\end{equation}

For $T$ a Young tableau, the \emph{row reading word of $T$}, denoted by $w(T)$, is the word obtained by reading the entries of $T$ along rows, left to right, from top to bottom. 

\begin{definition}
  The \emph{lowering operators}, denoted by $f_i$, on semistandard Young tableaux act by: $f_i(T)=0$ if $m_i(w(T)) \leqslant 0$; otherwise, letting $p$ be the smallest index such that $m_i(w(T),p) = m_i(w(T))$, $f_i(T)$ changes the entry in $T$ corresponding to $w_p$ to $i+1$.
  \label{def:young-lower}
\end{definition}

\begin{figure}[ht]
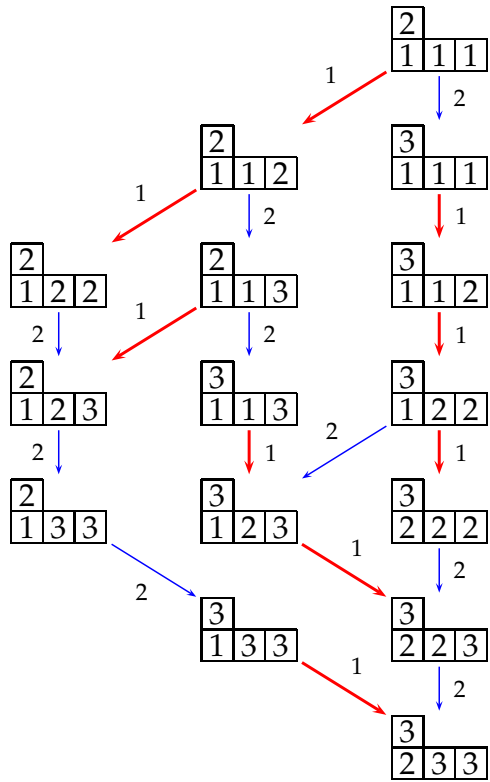

  \begin{center}
    \begin{displaymath}
      \begin{array}{c@{\hskip 3\cellsize}c@{\hskip 3\cellsize}c}
        & & \rnode{c1}{\tableau{2 \\ 1 & 1 & 1}} \\[5ex]
        & \rnode{b2}{\tableau{2 \\ 1 & 1 & 2}} & \rnode{c2}{\tableau{3 \\ 1 & 1 & 1}} \\[5ex]
        \rnode{a3}{\tableau{2 \\ 1 & 2 & 2}} & \rnode{b3}{\tableau{2 \\ 1 & 1 & 3}} & \rnode{c3}{\tableau{3 \\ 1 & 1 & 2}} \\[5ex]
        \rnode{a4}{\tableau{2 \\ 1 & 2 & 3}} & \rnode{b4}{\tableau{3 \\ 1 & 1 & 3}} & \rnode{c4}{\tableau{3 \\ 1 & 2 & 2}} \\[5ex]
        \rnode{a5}{\tableau{2 \\ 1 & 3 & 3}} & \rnode{b5}{\tableau{3 \\ 1 & 2 & 3}} & \rnode{c5}{\tableau{3 \\ 2 & 2 & 2}} \\[5ex]
        & \rnode{b6}{\tableau{3 \\ 1 & 3 & 3}} & \rnode{c6}{\tableau{3 \\ 2 & 2 & 3}} \\[5ex]
        & & \rnode{c7}{\tableau{3 \\ 2 & 3 & 3}}
      \end{array}
      \psset{nodesep=2pt,linewidth=.1ex}
      \everypsbox{\scriptstyle}
      \ncline[linewidth=.2ex,linecolor=red]{->} {c1}{b2} \nbput{1}
      \ncline[linecolor=blue]{->}  {c1}{c2} \naput{2}
      \ncline[linewidth=.2ex,linecolor=red]{->}   {b2}{a3} \nbput{1}
      \ncline[linecolor=blue]{->}  {b2}{b3} \naput{2}
      \ncline[linewidth=.2ex,linecolor=red]{->} {c2}{c3} \naput{1}
      \ncline[linecolor=blue]{->}  {a3}{a4} \nbput{2}
      \ncline[linewidth=.2ex,linecolor=red]{->}   {b3}{a4} \nbput{1}
      \ncline[linecolor=blue]{->}  {b3}{b4} \naput{2}
      \ncline[linewidth=.2ex,linecolor=red]{->} {c3}{c4} \naput{1}
      \ncline[linecolor=blue]{->}  {a4}{a5} \nbput{2}
      \ncline[linewidth=.2ex,linecolor=red]{->} {b4}{b5} \naput{1}
      \ncline[linecolor=blue]{->}  {c4}{b5} \nbput{2}
      \ncline[linewidth=.2ex,linecolor=red]{->} {c4}{c5} \naput{1}
      \ncline[linecolor=blue]{->}  {a5}{b6} \nbput{2}
      \ncline[linewidth=.2ex,linecolor=red]{<-}  {c6}{b5} \nbput{1} 
      \ncline[linecolor=blue]{->}  {c5}{c6} \naput{2}
      \ncline[linewidth=.2ex,linecolor=red]{<-}  {c7}{b6} \nbput{1} 
      \ncline[linecolor=blue]{->}  {c6}{c7} \naput{2}
    \end{displaymath}
    \caption{\label{fig:crystal-A}The crystal graph on $\SSYT(3,1)$ with entries $\{1,2,3\}$.}
  \end{center}
\end{figure}

Note that if $m_i(w) > 0$ and $w_p$ is the leftmost occurrence of this maximum, then $w_p = i$. In particular, this makes $f_i$ well-defined. For examples, see Figure~\ref{fig:crystal-A}. 

\begin{theorem}[\cite{KN94,Lit95}]
  The lowering operators $\{f_i\}_{1 \leq i < n}$ are well-defined maps $f_i : \SSYT_n(\lambda) \rightarrow \SSYT_n(\lambda) \cup \{0\}$ that define a crystal graph on semistandard Young tableaux.
  \label{thm:young-crystal}
\end{theorem}

Stembridge \cite{Ste03} gave a local characterization of the crystal graphs on semistandard Young tableaux in terms of six local axioms. This gives a powerful tool for establishing the symmetry and Schur positivity of a given polynomial by defining lowering operators on the generating objects and proving that they satisfy Stembridge's axioms.

Say that $\mathcal{X}$ has {\em degree $n$} if the largest color for an edge is $n-1$. If $x {\blue \stackrel{i}{\leftarrow}} y$ (resp. $x {\blue \stackrel{i}{\rightarrow}} z$), write $y=E_i x$ (resp. $z = F_i x$). The {\em $i$-string through $x$} is the maximal path
$$
F_i^{-d} x {\blue \stackrel{i}{\rightarrow}} \cdots {\blue \stackrel{i}{\rightarrow}} F_i^{-1} x {\blue \stackrel{i}{\rightarrow}} x
{\blue \stackrel{i}{\rightarrow}} F_i x {\blue \stackrel{i}{\rightarrow}} \cdots {\blue \stackrel{i}{\rightarrow}} F_i^{r} x .
$$
In this case we write $\delta_i(x) = d$ and $\varepsilon_i(x) = r$. Whenever $E_i, F_i$ is defined at $x$, we set
\begin{displaymath}
  \begin{array}{rclcrcl}
    \Delta_i \delta_j(x) & = & \delta_j( x) - \delta_j(E_ix), & &
    \nabla_i \delta_j(x) & = & \delta_j(F_ix) - \delta_j( x), \\
    \Delta_i \varepsilon_j(x) & = & \varepsilon_j(E_i x) - \varepsilon_j(x), & &
    \nabla_i \varepsilon_j(x) & = & \varepsilon_j(x) - \varepsilon_j(F_i x).
  \end{array}
\end{displaymath}

\begin{definition}[\cite{Ste03}]
  A directed, colored graph $\mathcal{X}$ is {\em regular} if the following hold:
  \begin{itemize}
    \addtolength{\itemsep}{-0.5\baselineskip}
  \item[(A1)] all monochromatic directed paths have finite length;

  \item[(A2)] for every vertex $x$, there is at most one edge $x {\blue \stackrel{i}{\leftarrow}} y$ and at most one edge $x {\blue \stackrel{i}{\rightarrow}} z$;

  \item[(A3)] assuming $E_i x$ is defined, $\Delta_i \delta_j(x) + \Delta_i \varepsilon_j(x) = \left\{
    \begin{array}{rl}
        2 & \;\mbox{if}\;\; j=i \\
        -1 & \;\mbox{if}\;\; j=i\pm 1 \\
        0 & \;\mbox{if}\;\; |i-j|\geq 2
    \end{array}
    \right.$;
    
  \item[(A4)] assuming $E_i x$ is defined, $\Delta_i \delta_j(x), \Delta_i \varepsilon_j(x) \leq 0$ for $j \neq i$;

  \item[(A5)] $\Delta_i \delta_j(x) = 0$ $\Rightarrow$ $E_iE_j x = E_jE_i x = y$ and $\nabla_j \varepsilon_i(y) = 0$; \\
    $\nabla_i \varepsilon_j(x) = 0$ $\Rightarrow$ $F_iF_j x = F_jF_i x = y$ and $\Delta_j \delta_i(y) = 0$;

  \item[(A6)] $\Delta_i \delta_j(x) = \Delta_j \delta_i(x) = -1$ $\Rightarrow$ $E_iE_{j}^{2}E_i x = E_jE_{i}^{2}E_j x = y$ and $\nabla_i \varepsilon_j(y) = \nabla_j \varepsilon_i(y)=-1$; \\
    $\nabla_i \varepsilon_j(x) = \nabla_j \varepsilon_i(x) = -1$ $\Rightarrow$ $F_iF_{j}^{2}F_i x = F_jF_{i}^{2}F_j x = y$ and $\Delta_i \delta_j(y) = \Delta_j \delta_i(y)=-1$.
  \end{itemize}
  \label{defn:regular}
\end{definition}

\begin{figure}[ht]
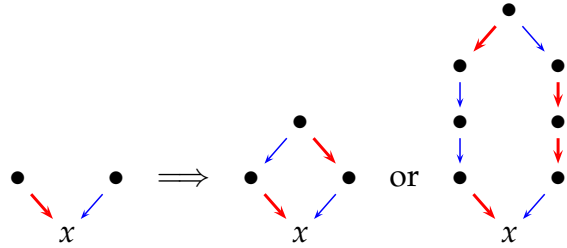

  \begin{center}
    \begin{displaymath}
      \begin{array}{c@{\hskip\cellsize}c@{\hskip\cellsize}c@{\hskip\cellsize}c@{\hskip\cellsize}c@{\hskip\cellsize}c@{\hskip\cellsize}c@{\hskip\cellsize}c@{\hskip\cellsize}c@{\hskip\cellsize}c@{\hskip\cellsize}c}  
        & & & & & & & & & \rnode{a11}{\bullet} & \\[0.5\cellsize]
        & & & & & & & & \rnode{b10}{\bullet} & & \rnode{b12}{\bullet} \\[0.5\cellsize]
        & & & & & \rnode{c6}{\bullet} & & & \rnode{c11l}{\bullet} & & \rnode{c11r}{\bullet} \\[0.5\cellsize] 
        \rnode{d1}{\bullet} & & \rnode{d3}{\bullet} & \rnode{d4}{\Longrightarrow} & \rnode{d5}{\bullet} & & \rnode{d7}{\bullet} & \rnode{d8}{\mbox{or}} & \rnode{d9}{\bullet} & & \rnode{d13}{\bullet} \\[0.5\cellsize]
        & \rnode{e2}{x} & & & & \rnode{e6}{x} & & & & \rnode{e11}{x} & 
      \end{array}
      \psset{linewidth=.1ex,nodesep=3pt}
      \everypsbox{\scriptstyle}
      \ncline[linecolor=red,linewidth=.2ex]{->} {d1}{e2}  
      \ncline[linecolor=blue]{->} {d3}{e2}  
      \ncline[linecolor=blue]{->} {c6}{d5}  
      \ncline[linecolor=red,linewidth=.2ex]{->} {c6}{d7}  
      \ncline[linecolor=red,linewidth=.2ex]{->} {d5}{e6}  
      \ncline[linecolor=blue]{->} {d7}{e6}  
      \ncline[linecolor=blue]{->} {a11}{b12}  
      \ncline[linecolor=red,linewidth=.2ex]{->} {a11}{b10}  
      \ncline[linecolor=red,linewidth=.2ex]{->} {b12}{c11r} 
      \ncline[linecolor=blue]{->} {b10}{c11l} 
      \ncline[linecolor=blue]{->} {c11l}{d9}  
      \ncline[linecolor=red,linewidth=.2ex]{->} {c11r}{d13} 
      \ncline[linecolor=red,linewidth=.2ex]{->} {d9}{e11}   
      \ncline[linecolor=blue]{->} {d13}{e11}  
    \end{displaymath}
    \caption{\label{fig:P5P6}An illustration of axioms A5 and A6, where $F_j {\blue \swarrow}$, $F_i {\red \searrow}$.}
  \end{center}
\end{figure}

\begin{theorem}[\cite{Ste03}]
  A graph is a regular graph if and only if every connected component is isomorphic to the crystal graph on $\SSYT_n(\lambda)$ for some $\lambda, n$.
\label{thm:structure-crystal}
\end{theorem}

%
\section{Crystal graph on semistandard shifted tableaux}
%
\label{sec:shifted}

A partition $\gamma$ is \emph{strict} if  $\gamma_1 > \gamma_2 > \cdots > \gamma_{\ell} > 0$. We identify a strict partition $\gamma$ with its \emph{shifted diagram}, the diagram with $\gamma_i$ cells in row $i$ shifted $\ell(\gamma)-i$ cells to the left. 

\begin{figure}[ht]
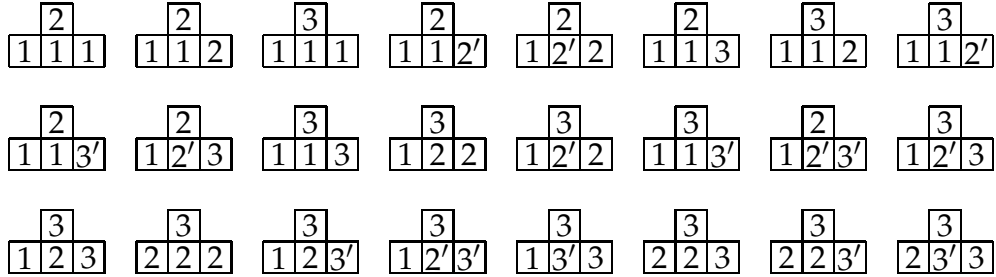

  \begin{center}
    \begin{displaymath}
      \begin{array}{c@{\hskip\cellsize}c@{\hskip\cellsize}c@{\hskip\cellsize}c@{\hskip\cellsize}c@{\hskip\cellsize}c@{\hskip\cellsize}c@{\hskip\cellsize}c}
        \tableau{ & 2 \\ 1 & 1 & 1} & \tableau{ & 2 \\ 1 & 1 & 2} & \tableau{ & 3 \\ 1 & 1 & 1} & \tableau{ & 2 \\ 1 & 1 & \st{2}} & \tableau{ & 2 \\ 1 & \st{2} & 2} & \tableau{ & 2 \\ 1 & 1 & 3} & \tableau{ & 3 \\ 1 & 1 & 2} & \tableau{ & 3 \\ 1 & 1 & \st{2}} \\[4ex]
        \tableau{ & 2 \\ 1 & 1 & \st{3}} & \tableau{ & 2 \\ 1 & \st{2} & 3} & \tableau{ & 3 \\ 1 & 1 & 3} & \tableau{ & 3 \\ 1 & 2 & 2} & \tableau{ & 3 \\ 1 & \st{2} & 2} & \tableau{ & 3 \\ 1 & 1 & \st{3}} & \tableau{ & 2 \\ 1 & \st{2} & \st{3}} & \tableau{ & 3 \\ 1 & \st{2} & 3} \\[4ex]
        \tableau{ & 3 \\ 1 & 2 & 3} & \tableau{ & 3 \\ 2 & 2 & 2} & \tableau{ & 3 \\ 1 & 2 & \st{3}} & \tableau{ & 3 \\ 1 & \st{2} & \st{3}} & \tableau{ & 3 \\ 1 & \st{3} & 3} & \tableau{ & 3 \\ 2 & 2 & 3} & \tableau{ & 3 \\ 2 & 2 & \st{3}} & \tableau{ & 3 \\ 2 & \st{3} & 3}
      \end{array}
    \end{displaymath}
    \caption{\label{fig:P31}The semistandard shifted tableaux for $(3,1)$ with entries $\{\st{1},1,\st{2},2,\st{3},3\}$.}
  \end{center}
\end{figure}

The \emph{semistandard shifted tableau} of shape $\gamma$, denoted by $\SSHT(\gamma)$, are fillings of the shifted diagram for $\gamma$ with marked or unmarked positive integers such that entries weakly increase along rows and columns according to the ordering $\st{1} < 1 < \st{2} < 2 < \cdots$, and for each $i$, at most one unmarked $i$ per column and at most one marked $\st{i}$ per row and none on the main diagonal. For example, see Figure~\ref{fig:P31}. 

\begin{definition}
  The Schur $P$-polynomial indexed by the strict partition $\gamma$ is given by
  \begin{equation}
    P_{\gamma} (x_1,\ldots,x_n) = \sum_{T \in \SSHT(\gamma)} x_1^{\wt(T)_1} \cdots x_n^{\wt(T)_n},
    \label{e:schur}
  \end{equation}
  where $\wt(T)_i$ is the number of occurrences of $i$ and $\st{i}$ in $T$.
  \label{def:schur}
\end{definition}
  
For example, from Figure~\ref{fig:P31}, we compute
\begin{eqnarray*}
  P_{(3,1)}(x_1,x_2,x_3) & = & x_1^3 x_2 + x_1^3 x_3 + 2 x_1^2 x_2^2 + 4 x_1^2 x_2 x_3 + 2 x_1^2 x_3^2 + x_1 x_2^3 \\
  & & + 4 x_1 x_2^2 x_3 + 4 x_1 x_2 x_3^2 + x_1 x_3^3 + x_2^3 x_3 + 2 x_2^2 x_3^2 + x_2 x_3^3.
\end{eqnarray*}

\begin{theorem}[\cite{Sag87,Wor84}]
  For $\gamma$ a strict partition, the coefficients $g_{\gamma,\lambda}$ defined by
    \begin{equation}
      P_{\gamma} (x_1,\ldots,x_n) = \sum_{\lambda} g_{\gamma,\lambda} s_{\lambda}(x_1,\ldots,x_n)
      \label{e:P_schur}
    \end{equation}
    are nonnegative integers. That is, Schur $P$-polynomials are Schur positive.
\label{thm:P_Schur}
\end{theorem}

We give a new proof of the Schur positivity of Schur $P$-polynomials by constructing a crystal graph on semistandard shifted tableaux. Hawkes, Paramonov, and Schilling recently gave a type A crystal for semistandard shifted tableaux \cite{HPS17} implicitly using mixed insertion. Our independent construction below coincides (though the proof of this is quite involved), but our proof is direct using Stembridge's axioms. 

\begin{definition}
  For $T$ a shifted tableau, the \emph{hook reading word of $T$}, denoted by $w(T)$, is the word obtained by reading the marked entries of $T$ up the $i$th column then the unmarked entries of $T$ along the $i$th row, left to right, for $i$ from $\max(\gamma_1,\ell(\gamma))$ to $1$.
  \label{def:hook-word}
\end{definition}

Note that the cells with entries $\st{i},i$ must form a ribbon, i.e. contain no $2\times 2$ block.

\begin{definition}
  The \emph{shifted lowering operators}, denoted by $\st{f}_i$, act on semistandard shifted tableaux by: $\st{f}_i(T)=0$ if $m_i(w(T)) \leqslant 0$; otherwise, letting $p$ be the smallest index such that $m_i(w(T),p) = m_i(w(T))$, letting $x$ denote the entry of $T$ corresponding to $w_p$, and letting $y$ be the entry north of $x$ and $z$ the entry east of $x$, we have
    \begin{enumerate}
    \item \begin{enumerate}
      \item if $x=i$ and $z=\st{i+1}$, then $\st{f}_i$ changes $x$ to $\st{i+1}$ and changes $z$ to $i+1$;
      \item else if $x=i$ and $y$ does not exist or $y > i+1$, then $\st{f}_i$ changes $x$ to $i+1$;
      \item else if $x=i$, then $\st{f}_i$ changes $x$ to $\st{i+1}$ and removes the marking (if any) from northwestern-most cell on the $(i+1)$-ribbon containing $y$;
      \end{enumerate}
      \begin{figure}[ht]
        \begin{displaymath}
          \begin{array}{ccc@{\hskip 4\cellsize}ccc@{\hskip 4\cellsize}ccc}
            \tableau{\\ y \\ 1 & \st{2} } & \stackrel{1(a)}{\mapsto} & \tableau{\\ y \\ \st{2} & 2 } &
            \tableau{\\ y \\ 1 & z }      & \stackrel{1(b)}{\mapsto} & \tableau{\\ y \\ 2 & z } &
            \tableau{\st{2} & 2 & 2 \\ & & \st{2} \\ & & 1 & z } & \stackrel{1(c)}{\mapsto} &
            \tableau{ 2 & 2 & 2 \\ & & \st{2} \\ & & \st{2} & z }  
          \end{array}
        \end{displaymath}
        \caption{\label{fig:shifted-lower-1}An illustration of the shifted lowering operators with $x=i=1$.}
      \end{figure}

    \item \begin{enumerate}
      \item if $x=\st{i}$ and $y = i$, then $\st{f}_i$ changes $x$ to $i$ and changes $y$ to $\st{i+1}$;
      \item else if $x=\st{i}$ and $z$ does not exist or $z > \st{i+1}$, then $\st{f}_i$ changes $x$ to $\st{i+1}$;
      \item else if $x=\st{i}$, then $\st{f}_i$ changes $x$ to $i$ and changes the first entry $i$ southwest along the $i$-ribbon containing $x$ that is not followed by $i$ or $\st{i+1}$ to $\st{i+1}$.
      \end{enumerate}
      \begin{figure}[ht]
        \begin{displaymath}
          \begin{array}{ccc@{\hskip 4\cellsize}ccc@{\hskip 4\cellsize}ccc}
            \tableau{ 1 \\ \st{1} & z } & \stackrel{2(a)}{\mapsto} & \tableau{ \st{2} \\ 1 & z } &
            \tableau{ y \\ \st{1} & z } & \stackrel{2(b)}{\mapsto} & \tableau{ y \\ \st{2} & z } &
            \tableau{ y \\ \st{1} & 1 & 1 & \st{2} \\ & & \st{1} & 1 } & \stackrel{2(c)}{\mapsto} &
            \tableau{ y \\ 1 & 1 & 1 & \st{2} \\ & & \st{1} & \st{2} } 
          \end{array}
        \end{displaymath}
        \caption{\label{fig:shifted-lower-2}An illustration of the shifted lowering operators with $x=\st{i}=\st{1}$.}
      \end{figure}
    \end{enumerate}
  \label{def:shifted-lower}
\end{definition}

\begin{figure}[ht]
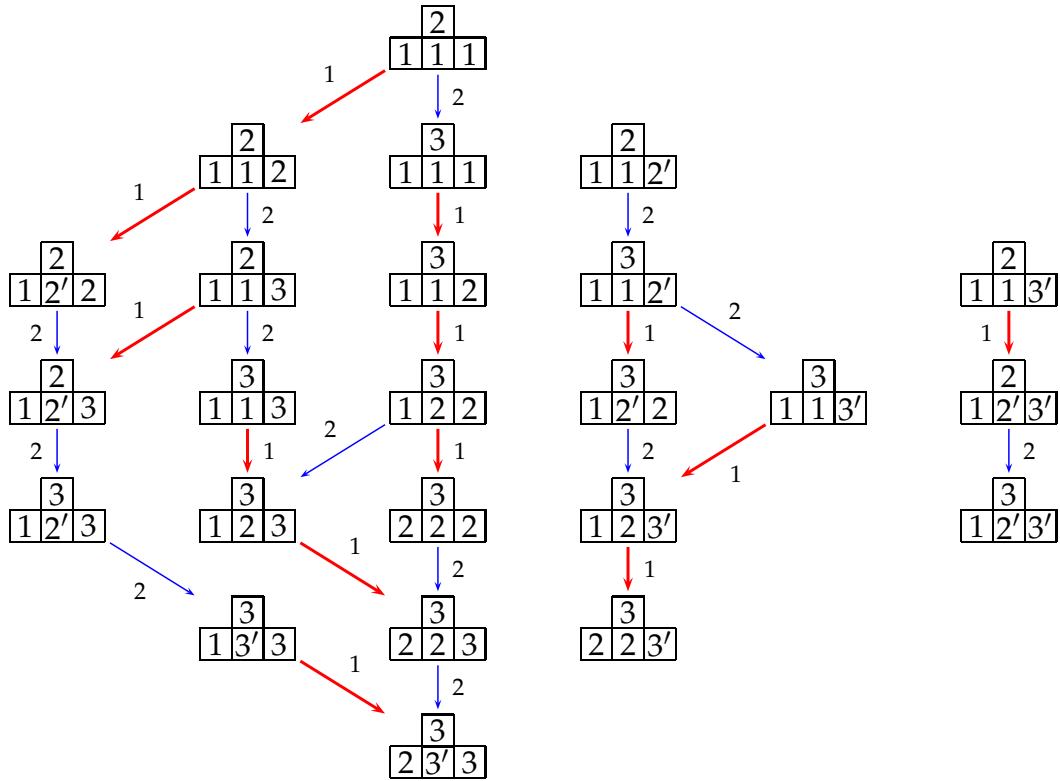

  \begin{center}
    \begin{displaymath}
      \begin{array}{c@{\hskip 3\cellsize}c@{\hskip 3\cellsize}c@{\hskip 3\cellsize}c@{\hskip 3\cellsize}c@{\hskip 3\cellsize}c}
        & & \rnode{c1}{\tableau{ & 2 \\ 1 & 1 & 1}} & & & \\[5ex]
        & \rnode{b2}{\tableau{ & 2 \\ 1 & 1 & 2}} & \rnode{c2}{\tableau{ & 3 \\ 1 & 1 & 1}} & \rnode{d2}{\tableau{ & 2 \\ 1 & 1 & \st{2}}} & & \\[5ex]
        \rnode{a3}{\tableau{ & 2 \\ 1 & \st{2} & 2}} & \rnode{b3}{\tableau{ & 2 \\ 1 & 1 & 3}} & \rnode{c3}{\tableau{ & 3 \\ 1 & 1 & 2}} & \rnode{d3}{\tableau{ & 3 \\ 1 & 1 & \st{2}}} & & \rnode{f3}{\tableau{ & 2 \\ 1 & 1 & \st{3}}}\\[5ex]
        \rnode{a4}{\tableau{ & 2 \\ 1 & \st{2} & 3}} & \rnode{b4}{\tableau{ & 3 \\ 1 & 1 & 3}} & \rnode{c4}{\tableau{ & 3 \\ 1 & 2 & 2}} & \rnode{d4}{\tableau{ & 3 \\ 1 & \st{2} & 2}} & \rnode{e4}{\tableau{ & 3 \\ 1 & 1 & \st{3}}} & \rnode{f4}{\tableau{ & 2 \\ 1 & \st{2} & \st{3}}}\\[5ex]
        \rnode{a5}{\tableau{ & 3 \\ 1 & \st{2} & 3}} & \rnode{b5}{\tableau{ & 3 \\ 1 & 2 & 3}} & \rnode{c5}{\tableau{ & 3 \\ 2 & 2 & 2}} & \rnode{d5}{\tableau{ & 3 \\ 1 & 2 & \st{3}}} & & \rnode{f5}{\tableau{ & 3 \\ 1 & \st{2} & \st{3}}}\\[5ex]
        & \rnode{b6}{\tableau{ & 3 \\ 1 & \st{3} & 3}} & \rnode{c6}{\tableau{ & 3 \\ 2 & 2 & 3}} & \rnode{d6}{\tableau{ & 3 \\ 2 & 2 & \st{3}}} & & \\[5ex]
        & & \rnode{c7}{\tableau{ & 3 \\ 2 & \st{3} & 3}} & & & 
      \end{array}
      \psset{nodesep=2pt,linewidth=.1ex}
      \everypsbox{\scriptstyle}
      \ncline[linewidth=.2ex,linecolor=red]{->} {c1}{b2} \nbput{1}
      \ncline[linecolor=blue]{->}  {c1}{c2} \naput{2}
      \ncline[linewidth=.2ex,linecolor=red]{->}   {b2}{a3} \nbput{1}
      \ncline[linecolor=blue]{->}  {b2}{b3} \naput{2}
      \ncline[linewidth=.2ex,linecolor=red]{->} {c2}{c3} \naput{1}
      \ncline[linecolor=blue]{->}  {d2}{d3} \naput{2}
      \ncline[linecolor=blue]{->}  {a3}{a4} \nbput{2}
      \ncline[linewidth=.2ex,linecolor=red]{->}   {b3}{a4} \nbput{1}
      \ncline[linecolor=blue]{->}  {b3}{b4} \naput{2}
      \ncline[linewidth=.2ex,linecolor=red]{->} {c3}{c4} \naput{1}
      \ncline[linewidth=.2ex,linecolor=red]{->} {d3}{d4} \naput{1}
      \ncline[linecolor=blue]{->}  {d3}{e4} \naput{2}
      \ncline[linewidth=.2ex,linecolor=red]{->}   {f3}{f4} \nbput{1}
      \ncline[linecolor=blue]{->}  {a4}{a5} \nbput{2}
      \ncline[linewidth=.2ex,linecolor=red]{->} {b4}{b5} \naput{1}
      \ncline[linecolor=blue]{->}  {c4}{b5} \nbput{2}
      \ncline[linewidth=.2ex,linecolor=red]{->} {c4}{c5} \naput{1}
      \ncline[linecolor=blue]{->}  {d4}{d5} \naput{2}
      \ncline[linewidth=.2ex,linecolor=red]{->} {e4}{d5} \naput{1}
      \ncline[linecolor=blue]{->}  {f4}{f5} \naput{2}
      \ncline[linecolor=blue]{->}  {a5}{b6} \nbput{2}
      \ncline[linewidth=.2ex,linecolor=red]{<-}  {c6}{b5} \nbput{1} 
      \ncline[linecolor=blue]{->}  {c5}{c6} \naput{2}
      \ncline[linewidth=.2ex,linecolor=red]{->} {d5}{d6} \naput{1}
      \ncline[linewidth=.2ex,linecolor=red]{<-}  {c7}{b6} \nbput{1} 
      \ncline[linecolor=blue]{->}  {c6}{c7} \naput{2}
    \end{displaymath}
    \caption{\label{fig:P31-A}The crystal graph on $\SSHT(3,1)$ with entries $\{\st{1},1,\st{2},2,\st{3},3\}$.}
  \end{center}
\end{figure}

The rules for the shifted lowering operators are illustrated case by case in Figures~\ref{fig:shifted-lower-1} and \ref{fig:shifted-lower-2}, with examples on semistandard shifted tableaux shown in Figure~\ref{fig:P31-A}. 

\begin{theorem}
  The shifted lowering operators $\{f_i\}_{1 \leq i < n}$ are well-defined maps $\st{f}_i : \SSHT_n(\gamma) \rightarrow \SSHT_n(\gamma) \cup \{0\}$. Moreover, $\st{f}_i$ is injective on $\{ T \in \SSHT_n(\gamma) \mid \st{f}_i(T) \neq 0 \}$.
  \label{thm:shifted-well}
\end{theorem}

Using Stembridge's local characterization of the crystal graphs on Young tableaux \cite{Ste03}, we establish that the shifted lowering operators define a crystal on shifted tableaux.

\begin{theorem}
  The shifted lowering operators satisfy Stembridge's axioms. In particular, they define a crystal graph on semistandard shifted tableaux.
  \label{thm:shifted-crystal}
\end{theorem}

From this, we give a new proof of Theorem~\ref{thm:P_Schur}. Note that this characterization of the Schur coefficients of a Schur $P$-polynomial is simpler than Sagan's shifted insertion rule \cite{Sag87} and more explicit than Assaf's dual equivalence characterization \cite{Ass18}.

\begin{corollary}
  For $\gamma$ a strict partition, let $\Yam(\gamma)$ denote the set of semistandard shifted tableaux $T$ of shape $\gamma$ such that $m_i(T)=0$ for all $i$. Then
  \begin{equation}
    P_{\gamma}(x_1,\ldots,x_n) = \sum_{T \in \Yam(\gamma)} s_{\wt(T)}(x_1,\ldots,x_n).
  \end{equation}
  In particular, the Schur $P$-polynomial is Schur positive.
  \label{cor:P2schur}
\end{corollary}

For example, from Figure~\ref{fig:P31-A}, we have
\begin{displaymath}
  P_{(3,1)}(x_1,x_2,x_3) = s_{(3,1)}(x_1,x_2,x_3) + s_{(2,2)}(x_1,x_2,x_3) + s_{(2,1,1)}(x_1,x_2,x_3).
\end{displaymath}

%
\section{Queer crystal graph on semistandard shifted tableaux}
%
\label{sec:queer}

Since semistandard shifted tableaux describe tensor representations of the queer Lie superalgebra, the connected queer crystal that combinatorializes these representations becomes a powerful tool in proving Schur $P$-positivity. We may extend our crystal graph to a queer crystal by augmenting our lowering operators as follows.

\begin{definition}
  The \emph{queer lowering operator}, denoted by $\st{f}_0$, acts on semistandard shifted tableaux by: $\st{f}_0(T)=0$ if $T$ has no entry equal to $1$ or has an entry equal to $\st{2}$; otherwise $\st{f}_0(T)$ changes the rightmost $1$ in $T$ to $\st{2}$ if not on the main diagonal or to $2$ otherwise.
  \label{def:queer-lower}
\end{definition}

Figure~\ref{fig:P31-B} demonstrates the queer lowering operator on semistandard shifted tableaux.

\begin{figure}[ht]
  \begin{center}
    \begin{displaymath}
      \begin{array}{c@{\hskip 3\cellsize}c@{\hskip 3\cellsize}c@{\hskip 3\cellsize}c@{\hskip 3\cellsize}c@{\hskip 3\cellsize}c}
        & & \rnode{c1}{\tableau{ & 2 \\ 1 & 1 & 1}} & & & \\[5ex]
        & \rnode{b2}{\tableau{ & 2 \\ 1 & 1 & 2}} & \rnode{c2}{\tableau{ & 3 \\ 1 & 1 & 1}} & \rnode{d2}{\tableau{ & 2 \\ 1 & 1 & \st{2}}} & & \\[5ex]
        \rnode{a3}{\tableau{ & 2 \\ 1 & \st{2} & 2}} & \rnode{b3}{\tableau{ & 2 \\ 1 & 1 & 3}} & \rnode{c3}{\tableau{ & 3 \\ 1 & 1 & 2}} & \rnode{d3}{\tableau{ & 3 \\ 1 & 1 & \st{2}}} & & \rnode{f3}{\tableau{ & 2 \\ 1 & 1 & \st{3}}}\\[5ex]
        \rnode{a4}{\tableau{ & 2 \\ 1 & \st{2} & 3}} & \rnode{b4}{\tableau{ & 3 \\ 1 & 1 & 3}} & \rnode{c4}{\tableau{ & 3 \\ 1 & 2 & 2}} & \rnode{d4}{\tableau{ & 3 \\ 1 & \st{2} & 2}} & \rnode{e4}{\tableau{ & 3 \\ 1 & 1 & \st{3}}} & \rnode{f4}{\tableau{ & 2 \\ 1 & \st{2} & \st{3}}}\\[5ex]
        \rnode{a5}{\tableau{ & 3 \\ 1 & \st{2} & 3}} & \rnode{b5}{\tableau{ & 3 \\ 1 & 2 & 3}} & \rnode{c5}{\tableau{ & 3 \\ 2 & 2 & 2}} & \rnode{d5}{\tableau{ & 3 \\ 1 & 2 & \st{3}}} & & \rnode{f5}{\tableau{ & 3 \\ 1 & \st{2} & \st{3}}}\\[5ex]
        & \rnode{b6}{\tableau{ & 3 \\ 1 & \st{3} & 3}} & \rnode{c6}{\tableau{ & 3 \\ 2 & 2 & 3}} & \rnode{d6}{\tableau{ & 3 \\ 2 & 2 & \st{3}}} & & \\[5ex]
        & & \rnode{c7}{\tableau{ & 3 \\ 2 & \st{3} & 3}} & & & 
      \end{array}
      \psset{nodesep=2pt,linewidth=.1ex}
      \everypsbox{\scriptstyle}
      \ncline[linewidth=.2ex,linecolor=red]{->} {c1}{b2} \nbput{1}
      \ncline[linecolor=blue]{->}  {c1}{c2} \naput{2}
      \ncline[linewidth=.3ex,linecolor=Green]{->}   {c1}{d2} \naput{0}
      \ncline[offset=2pt,linewidth=.3ex,linecolor=Green]{->}   {b2}{a3} \nbput{1}
      \ncline[offset=2pt,linewidth=.2ex,linecolor=red]{<-} {a3}{b2} \nbput{0}
      \ncline[linecolor=blue]{->}  {b2}{b3} \naput{2}
      \ncline[linewidth=.2ex,linecolor=red]{->} {c2}{c3} \naput{1}
      \ncline[linewidth=.3ex,linecolor=Green]{->}   {c2}{d3} \naput{0}
      \ncline[linecolor=blue]{->}  {d2}{d3} \naput{2}
      \ncline[linecolor=blue]{->}  {a3}{a4} \nbput{2}
      \ncline[offset=2pt,linewidth=.3ex,linecolor=Green]{->}   {b3}{a4} \nbput{1}
      \ncline[offset=2pt,linewidth=.2ex,linecolor=red]{<-} {a4}{b3} \nbput{0} 
      \ncline[linecolor=blue]{->}  {b3}{b4} \naput{2}
      \ncline[linewidth=.2ex,linecolor=red]{->} {c3}{c4} \naput{1}
      \ncline[linewidth=.3ex,linecolor=Green]{->}   {c3}{d4} \naput{0}
      \ncline[linewidth=.2ex,linecolor=red]{->} {d3}{d4} \naput{1}
      \ncline[linecolor=blue]{->}  {d3}{e4} \naput{2}
      \ncline[offset=2pt,linewidth=.3ex,linecolor=Green]{->}   {f3}{f4} \nbput{1}
      \ncline[offset=2pt,linewidth=.2ex,linecolor=red]{<-} {f4}{f3} \nbput{0} 
      \ncline[linecolor=blue]{->}  {a4}{a5} \nbput{2}
      \ncline[linewidth=.3ex,linecolor=Green]{->}   {b4}{a5} \naput{0}
      \ncline[linewidth=.2ex,linecolor=red]{->} {b4}{b5} \naput{1}
      \ncline[linecolor=blue]{->}  {c4}{b5} \nbput{2}
      \ncline[offset=2pt,linewidth=.3ex,linewidth=.3ex,linecolor=Green]{->} {c4}{c5} \nbput{1}
      \ncline[offset=2pt,linewidth=.2ex,linecolor=red]{<-} {c5}{c4} \nbput{0}
      \ncline[linecolor=blue]{->}  {d4}{d5} \naput{2}
      \ncline[linewidth=.2ex,linecolor=red]{->} {e4}{d5} \naput{1}
      \ncline[linewidth=.3ex,linecolor=Green]{->}   {e4}{f5} \naput{0}
      \ncline[linecolor=blue]{->}  {f4}{f5} \naput{2}
      \ncline[linecolor=blue]{->}  {a5}{b6} \nbput{2}
      \ncline[offset=2pt,linewidth=.2ex,linecolor=red]{->}{b5}{c6} \nbput{0}
      \ncline[offset=2pt,linewidth=.3ex,linecolor=Green]{<-}  {c6}{b5} \nbput{1} 
      \ncline[linecolor=blue]{->}  {c5}{c6} \naput{2}
      \ncline[offset=2pt,linewidth=.3ex,linecolor=Green]{->} {d5}{d6} \nbput{1}
      \ncline[offset=2pt,linewidth=.2ex,linecolor=red]{<-} {d6}{d5} \nbput{0}
      \ncline[offset=2pt,linewidth=.2ex,linecolor=red]{->}{b6}{c7} \nbput{0}
      \ncline[offset=2pt,linewidth=.3ex,linecolor=Green]{<-}  {c7}{b6} \nbput{1} 
      \ncline[linecolor=blue]{->}  {c6}{c7} \naput{2}
    \end{displaymath}
    \caption{\label{fig:P31-B}The queer crystal graph on $\SSHT(3,1)$ with entries $\{\st{1},1,\st{2},2,\st{3},3\}$.}
  \end{center}
\end{figure}

Using the definition of an abstract queer crystal from \cite{GJKKK14}, we have the following.

\begin{theorem}
  The queer and shifted lowering operators $\{f_i\}_{0 \leq i < n}$ define an abstract queer crystal graph on semistandard Young tableaux.
  \label{thm:queer-crystal}
\end{theorem}

Stembridge \cite{Ste89} expanded on Sagan's shifted insertion \cite{Sag87} in his study of projective representations of the symmetric group, and ultimately established that the product of Schur $P$ functions expands positively in the Schur $P$ basis. More recently, Cho \cite{Cho13} built on work of Serrano \cite{Ser10} to give another proof of positivity, and Assaf \cite{Ass18} gave another proof using the machinery of combinatorial graphs.

Using the combinatorial rule for tensoring queer crystals together with the characterization of highest weights given in \cite{GJKKK14}, we have the following explicit expansion.

\begin{corollary}
  For $\gamma,\delta$ strict partitions, the coefficients $f_{\gamma,\delta}^{\varepsilon}$ defined by
    \begin{equation}
      P_{\gamma} (x_1,\ldots,x_n) P_{\delta} (x_1,\ldots,x_n) = \sum_{\varepsilon} f_{\gamma,\delta}^{\varepsilon} P_{\varepsilon}(x_1,\ldots,x_n),
      \label{e:P_product}
    \end{equation}
    where $f_{\gamma,\delta}^{\varepsilon}$ is the number highest weight tableaux $S \otimes T$ of highest weight $\varepsilon$ in the tensor product of the queer crystal graphs $B(\gamma) \otimes B(\delta)$. In particular, $f_{\gamma,\delta}^{\varepsilon}$ is a nonnegative integer. That is, products of Schur $P$ polynomials are Schur $P$ positive.
  \label{cor:LRR}
\end{corollary}

A queer crystal is a regular crystal together with additional queer edges. Following Stembridge \cite{Ste03}, we give a local characterization of the queer crystal graphs on semistandard shifted tableaux by augmenting Stembridge's axioms with the following.

\begin{definition}
  A directed, colored graph $\mathcal{X}$ is {\em queer regular} if the following hold:
  \begin{itemize}
    \addtolength{\itemsep}{-0.5\baselineskip}
  \item[(B0)] the subgraph $\mathcal{X}_+$ generated by all edges with non-zero labels is regular;
    
  \item[(B1)] all $0$ paths have length $1$;
    \\ 
    if $\delta_1(x)+\varepsilon_1(x)+\varepsilon_2(x)>1$ then $\delta_0(x)+\varepsilon_0(x)=1$;
    
  \item[(B2)] for every vertex $x$, there is at most one edge $x \stackrel{0}{\longleftarrow} y$ and at most one edge $x \stackrel{0}{\longrightarrow} z$;
    
  \item[(B3)] assuming $E_0 x$ is defined, $\Delta_0 \delta_i(x) + \Delta_0 \varepsilon_i(x) = \left\{
    \begin{array}{rl}
      2 & \;\mbox{if}\;\; i \leq 1 \\
      -1 & \;\mbox{if}\;\; i=2 \\
        0 & \;\mbox{if}\;\; i\geq 3
    \end{array}
    \right.$;
    
  \item[(B4)] assuming $E_0 x$ is defined,   $ \begin{array}{rl}
      \Delta_0 \delta_i(x)>0, \Delta_0 \varepsilon_i(x) \geq 0 & \;\mbox{if}\;\; i=1 \\
    \Delta_0 \delta_i(x)\leq 0 , \Delta_0 \varepsilon_i(x) \leq 0 & \;\mbox{if}\;\; i = 2 \\
    \Delta_0 \delta_i(x)=0, \Delta_0 \varepsilon_i(x) = 0 & \;\mbox{if}\;\; i\geq 3 
  \end{array} $
    
  \item[(B5)]For $i\geq 2$, $E_i x=E_0 y=z \Rightarrow F_i F_0 z=F_0 F_i z$; \\
     For $i=1$ or $i \geq 3$, $F_i x=F_0 y=z$ and $x\neq y \Rightarrow E_i E_0 z=E_0 E_i z$ ;
    
  \item[(B6)] assuming $E_0 x$ is defined,$\begin{array}{ll}
    \Delta_0 \delta_1(x)=1 \Rightarrow \varepsilon_1(x)=0 \text{ and } E_0x=E_1x\\
    \Delta_0 \varepsilon_2(x)=0 \Leftrightarrow \varepsilon_2(x)=0
  \end{array}$. 
  \end{itemize}
  \label{def:characterization}
\end{definition}

\begin{figure}[ht]
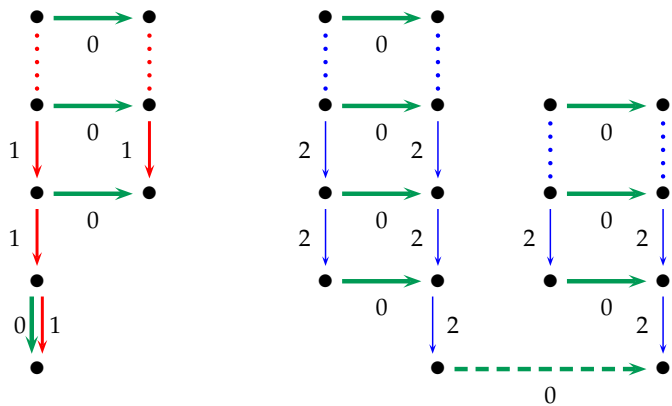

  \begin{displaymath}
    \begin{array}{c@{\hskip 3\cellsize}c@{\hskip 5\cellsize}c@{\hskip 3\cellsize}c@{\hskip 3\cellsize}c@{\hskip 3\cellsize}c}
      \rnode{A1}{\bullet} & \rnode{B1}{\bullet} & \rnode{a1}{\bullet} & \rnode{b1}{\bullet} &                     &                     \\[1.5\cellsize]
      \rnode{A2}{\bullet} & \rnode{B2}{\bullet} & \rnode{a2}{\bullet} & \rnode{b2}{\bullet} & \rnode{c2}{\bullet} & \rnode{d2}{\bullet} \\[1.5\cellsize]
      \rnode{A4}{\bullet} & \rnode{B4}{\bullet} & \rnode{a4}{\bullet} & \rnode{b4}{\bullet} & \rnode{c4}{\bullet} & \rnode{d4}{\bullet} \\[1.5\cellsize]
      \rnode{A5}{\bullet} &                     & \rnode{a5}{\bullet} & \rnode{b5}{\bullet} & \rnode{c5}{\bullet} & \rnode{d5}{\bullet} \\[1.5\cellsize]
      \rnode{A6}{\bullet} &                     &                     & \rnode{b6}{\bullet} &                     & \rnode{d6}{\bullet}               
    \end{array} 
    \psset{linewidth=.1ex,nodesep=3pt}
    \everypsbox{\scriptstyle}
    \ncline[linecolor=red,linestyle=dotted,linewidth=.3ex] {A1}{A2}
    \ncline[linecolor=red,linestyle=dotted,linewidth=.3ex] {B1}{B2}
    \ncline[linewidth=.2ex,linecolor=red]{->}              {A2}{A4} \nbput{1}
    \ncline[linewidth=.2ex,linecolor=red]{->}              {A4}{A5} \nbput{1}
    \ncline[linewidth=.2ex,linecolor=red]{->}              {B2}{B4} \nbput{1}
    \ncline[linewidth=.3ex,linecolor=ForestGreen]{->}              {A1}{B1} \nbput{0}
    \ncline[linewidth=.3ex,linecolor=ForestGreen]{->}              {A2}{B2} \nbput{0}
    \ncline[linewidth=.3ex,linecolor=ForestGreen]{->}              {A4}{B4} \nbput{0}
    \ncline[offset=2pt,linewidth=.2ex,linecolor=red]{->}{A5}{A6} \nbput{0}
    \ncline[offset=2pt,linewidth=.3ex,linecolor=ForestGreen]{<-}  {A6}{A5} \nbput{1}      
    \ncline[linecolor=blue,linestyle=dotted,linewidth=.3ex] {a1}{a2}
    \ncline[linecolor=blue,linestyle=dotted,linewidth=.3ex] {b1}{b2}
    \ncline[linecolor=blue,linestyle=dotted,linewidth=.3ex] {c2}{c4}
    \ncline[linecolor=blue,linestyle=dotted,linewidth=.3ex] {d2}{d4}
    \ncline[linecolor=blue]{->}              {a2}{a4} \nbput{2}
    \ncline[linecolor=blue]{->}              {a4}{a5} \nbput{2}
    \ncline[linecolor=blue]{->}              {c4}{c5} \nbput{2}
    \ncline[linecolor=blue]{->}              {d5}{d6} \nbput{2}
    \ncline[linecolor=blue]{->}              {d4}{d5} \nbput{2}
    \ncline[linecolor=blue]{->}              {b2}{b4} \nbput{2}
    \ncline[linecolor=blue]{->}              {b4}{b5} \nbput{2}
    \ncline[linewidth=.3ex,linecolor=ForestGreen]{->}              {a5}{b5} \nbput{0}
    \ncline[linewidth=.3ex,linecolor=ForestGreen]{->}              {a1}{b1} \nbput{0}
    \ncline[linewidth=.3ex,linecolor=ForestGreen]{->}              {a2}{b2} \nbput{0}
    \ncline[linewidth=.3ex,linecolor=ForestGreen]{->}              {a4}{b4} \nbput{0}
    \ncline[linewidth=.3ex,linecolor=ForestGreen]{->}              {c5}{d5} \nbput{0}
    \ncline[linewidth=.3ex,linecolor=ForestGreen]{->}              {c2}{d2} \nbput{0}
    \ncline[linewidth=.3ex,linecolor=ForestGreen]{->}              {c4}{d4} \nbput{0}
    \ncline[linewidth=.3ex,linecolor=ForestGreen,linestyle=dashed]{->}              {b6}{d6} \nbput{0}
    \ncline[offset=2pt,linecolor=blue]{<-}  {b6}{b5} \nbput{2}      
  \end{displaymath}
  \caption{\label{fig:queer-char}Illustration of the queer axioms relating queer and regular edges.}
\end{figure}

\begin{theorem}
  The queer crystal on semistandard shifted tableaux is a queer regular graph.
  \label{thm:characterization}
\end{theorem}

We anticipate that we will soon prove that the converse holds as well, namely that every connected component of a queer regular graph is isomorphic to the queer crystal on $\SSHT_n(\gamma)$ for some strict partition $\gamma$ and some positive integer $n$.

%
\section{Rectification of shifted tableaux}
%
\label{sec:rectify}

In \cite{GJKKK14}, the authors construct an explicit queer crystal graph on Serrano's semistandard decomposition tableaux \cite{Ser10}. They chose this less common indexing set claiming that ``the set of shifted semistandard Young tableaux of fixed shape does not have a natural crystal structure.'' We argue that the shifted crystal operators of Definition~\ref{def:shifted-lower} are indeed natural by defining a rectification map from semistandard shifted tableaux to semistandard Young tableaux that commutes with the corresponding crystal operators.

\begin{figure}[ht]
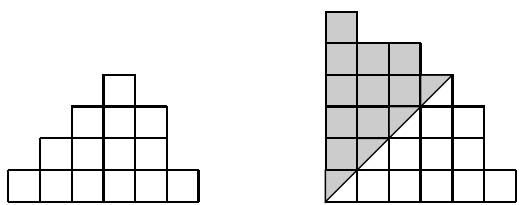

  \begin{displaymath}
    \tableau{\\ \\  & & & \ \\ &  & \ & \ & \ \\ & \ & \ & \ & \ \\ \ & \ & \ & \ & \ & \ }
    \hspace{4em}
    \tableau{ \cb \\ \cb & \cb & \cb \\ \cb & \cb & \cb & \tb \\ \cb & \cb & \tb & \ & \ \\ \cb & \tb & \ & \ & \ \\ \tb & \ & \ & \ & \ & \ }
  \end{displaymath}
  \caption{\label{fig:symmetric}The shifted (left) and symmetric (right) diagrams for $(6,4,3,1)$.}
\end{figure}

A partition $\nu$ is \emph{self-conjugate} if its diagram is symmetric about the main diagonal. For $\gamma$ a strict partition, the \emph{symmetric diagram} for $\gamma$, denoted by $\sym(\gamma)$, is the shifted diagram for $\gamma$ union with its transpose. For example, see Figure~\ref{fig:symmetric}. 

\begin{figure}[ht]
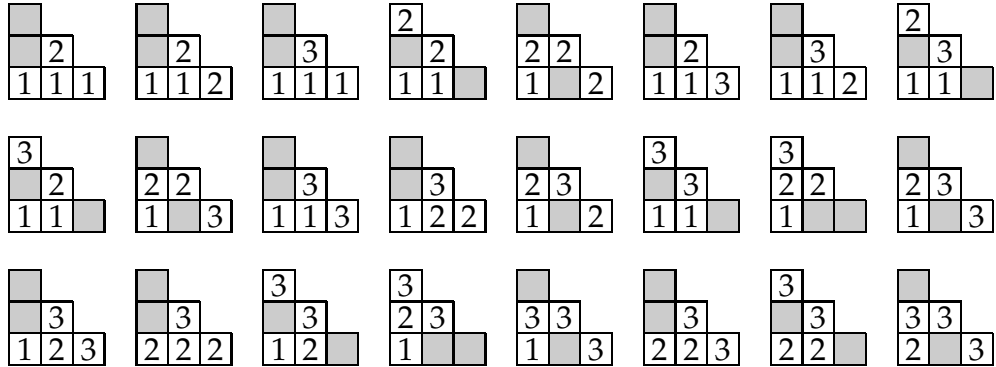

  \begin{center}
    \begin{displaymath}
      \begin{array}{c@{\hskip\cellsize}c@{\hskip\cellsize}c@{\hskip\cellsize}c@{\hskip\cellsize}c@{\hskip\cellsize}c@{\hskip\cellsize}c@{\hskip\cellsize}c}
        \tableau{ \cb \\ \cb & 2 \\ 1 & 1 & 1} &
        \tableau{ \cb \\ \cb & 2 \\ 1 & 1 & 2} &
        \tableau{ \cb \\ \cb & 3 \\ 1 & 1 & 1} &
        \tableau{ 2 \\ \cb & 2 \\ 1 & 1 & \cb} &
        \tableau{ \cb \\ 2 & 2 \\ 1 & \cb & 2} &
        \tableau{ \cb \\ \cb & 2 \\ 1 & 1 & 3} &
        \tableau{ \cb \\ \cb & 3 \\ 1 & 1 & 2} &
        \tableau{ 2 \\ \cb & 3 \\ 1 & 1 & \cb} \\[6ex]
        \tableau{ 3 \\ \cb & 2 \\ 1 & 1 & \cb} &
        \tableau{ \cb \\ 2 & 2 \\ 1 & \cb & 3} &
        \tableau{ \cb \\ \cb & 3 \\ 1 & 1 & 3} &
        \tableau{ \cb \\ \cb & 3 \\ 1 & 2 & 2} &
        \tableau{ \cb \\ 2 & 3 \\ 1 & \cb & 2} &
        \tableau{ 3 \\ \cb & 3 \\ 1 & 1 & \cb} &
        \tableau{ 3 \\ 2 & 2 \\ 1 & \cb & \cb} &
        \tableau{ \cb \\ 2 & 3 \\ 1 & \cb & 3} \\[6ex]
        \tableau{ \cb \\ \cb & 3 \\ 1 & 2 & 3} &
        \tableau{ \cb \\ \cb & 3 \\ 2 & 2 & 2} &
        \tableau{ 3 \\ \cb & 3 \\ 1 & 2 & \cb} &
        \tableau{ 3 \\ 2 & 3 \\ 1 & \cb & \cb} &
        \tableau{ \cb \\ 3 & 3 \\ 1 & \cb & 3} &
        \tableau{ \cb \\ \cb & 3 \\ 2 & 2 & 3} &
        \tableau{ 3 \\ \cb & 3 \\ 2 & 2 & \cb} &
        \tableau{ \cb \\ 3 & 3 \\ 2 & \cb & 3}
      \end{array}
    \end{displaymath}
    \caption{\label{fig:S31}The semistandard symmetric tableaux for $\sym(3,1)$ with entries $\{1,2,3\}$.}
  \end{center}
\end{figure}

The \emph{semistandard symmetric tableaux} of self-conjugate shape $\nu$, denoted by $\SSST(\nu)$, are partial fillings of the diagram for $\nu$ such that each diagonal cell is filled, for each non-diagonal cell $(j,k) \in \nu$ with $j \neq k$, exactly one of $(j,k)$ and $(k,j)$ is filled, and when entries in above-diagonal cells $(j,k)$ with $j < k$ are reflected to position $(k,j)$, the resulting filling satisfies the row and column conditions for Young tableaux.
For example, see Figure~\ref{fig:S31}. 

\begin{proposition}
  The map on semistandard shifted tableaux that reflects each marked entry in position $(j,k)$ to be an unmarked entry in position $(k,j)$ of the corresponding symmetric shape is a weight-preserving bijection $\SSHT(\gamma) \stackrel{\sim}{\rightarrow} \SSST(\sym(\gamma))$.
  \label{prop:sym}
\end{proposition}

For example, compare Figures~\ref{fig:P31} and \ref{fig:S31}, respectively. Notice that the hook reading word of a semistandard shifted tableau precisely corresponds to the usual row reading word of the corresponding semistandard symmetric tableau. 

\begin{definition}
  The \emph{rectification map}, denoted by $\rect(T)$, acts on partial fillings of partition shape with nonempty diagonal cells by: letting $d$ be the maximum row plus column coordinates of empty cells of $T$, letting $x$ be an empty cell at $(c,r)$ with $r+c=d$ such that $r$ is maximal (resp. minimal) among such empty cells if $r>c$ (resp. $r<c$), we have
  \begin{enumerate}
    \addtolength{\itemsep}{-0.6\baselineskip}
  \item if there is no cell above and no cell to the right, then delete $x$;
  \item if the cell above is filled with $i$ and either there is no cell to the right or it is filled with $j \geq i$, then fill $x$ with $i$ and empty the cell above;
  \item if the cell to the right is filled with $i$ and either there is no cell above or it is filled with $j > i$, then fill $x$ with $i$ and empty the cell to the right;
  \end{enumerate}
  then continue to rectify the newly vacated cell until it is deleted from the shape.
\label{def:rectify}
\end{definition} 


To see an example of rectification together with an illustration of the following theorem, the images of the semistandard shifted tableaux of shape $(3,1)$ on the leftmost component of the crystal in Figure~\ref{fig:P31-A} under the map to semistandard symmetric tableaux rectify to the semistandard Young tableaux of shape $(3,1)$ in Figure~\ref{fig:crystal-A}, respectively. 

\begin{theorem}
  Rectification induces a well-defined, weight-preserving bijection
  \[ \rect:\SSHT(\gamma) \stackrel{\sim}{\rightarrow} \bigsqcup_{T \in \Yam(\gamma)} \SSYT(\wt(T)). \]
  Moreover, for $T \in \SSHT(\gamma)$ and $i>0$, we have $\rect(\st{f}_i(T)) = f_i(\rect(T))$.
  \label{thm:rectify}
\end{theorem}


\printbibliography

\end{document}